\documentclass[12pt,a4paper]{amsart}%
\usepackage[utf8]{inputenc}
\usepackage{hyperref}
\usepackage{amsmath}
\usepackage{amsfonts}
\usepackage{amssymb}
\usepackage{xcolor}
\usepackage{graphicx}%
\providecommand{\U}[1]{\protect\rule{.1in}{.1in}}
\newtheorem{theorem}{Theorem}

\newtheorem{corollary}[theorem]{Corollary}

\newtheorem{lemma}[theorem]{Lemma}

\newtheorem{proposition}[theorem]{Proposition}
\newtheorem{remark}[theorem]{Remark}

\graphicspath{{Slike/}}

\begin{document}

\title{Upper bound $2^{n}-1$ for chromatic number of orbit polytope of a small cover is attained for every $n\geq4$}
\author{Djordje Barali\'{c}, Riste \v{S}krekovski}
\address{ \scriptsize{Mathematical Institute SANU\\ Knez Mihailova 36, p.f.\ 367\\
        11001 Belgrade, Serbia}}
\email{djbaralic@mi.sanu.ac.rs}
\address{ \scriptsize{Faculty of Information Studies\\ Ljubljanska cesta 31a\\
        8000 Novo Mesto, Slovenia}}
\address{ \scriptsize{Faculty of Mathematics and Physics, University of Ljubljana\\ Jadranska 19\\
        1000 Ljubljana, Slovenia}}
\email{skrekovski@gmail.com}

\date{}
\maketitle

\begin{abstract}
Davis and Januszkiewicz proved that the facets of a simple $n$-polytope
$P$ which is the orbit space of a small cover can be properly coloured by
$2^{n}-1$ colours. The bound is sharp for $n=2$, it is not sharp for $n=3$
by the Four Color Theorem, and it was shown to be sharp for $n=4$ by the
first author. We show that it is sharp for every $n\geq4$. The polytope is
obtained from the dual of the cyclic polytope $C^{n}(2^{n}-1)$ by
truncations of faces of codimension at least $3$. The resolution works over an
arbitrary field and gives the sharp bound $(q^{n}-1)/(q-1)$ over the field
with $q$ elements. More generally, every simple $n$-polytope whose facets can
be properly coloured by $k\leq2^{n}-1$ colours becomes, after truncations of
faces of codimension at least $3$, the orbit polytope of a small cover with
the same chromatic number $k$. We also construct,
for every $n\geq4$, a simple $n$-polytope with $\chi=2^{n}-1$ which is the
orbit space of both a small cover and a quasitoric manifold, and in the
oriented case we obtain the sharp bound $2^{n-1}$ for every $n\geq4$.
\end{abstract}



\section{Introduction}

Small covers and quasitoric manifolds were introduced by Davis and
Januszkiewicz \cite{DJ91} as topological generalizations of real and
complex smooth projective toric varieties. The orbit space of a small
cover is a simple polytope $P^{n}$. The manifold is determined by a
characteristic map, which assigns to each facet of $P^{n}$ a non-zero
vector of $\mathbb{Z}_{2}^{n}$ so that the vectors of the $n$ facets
meeting at a vertex form a basis of $\mathbb{Z}_{2}^{n}$, see \cite{BP15}.

A characteristic map is in particular a proper colouring of facets, where
no two facets sharing a vertex have the same colour. Hence the chromatic
number $\chi(P^{n})$ of a polytope admitting a characteristic map
satisfies \cite{DJ91} (here $\chi(P)$ is the chromatic number of the orbit
polytope $P$)
\[
n\leq\chi(P^{n})\leq2^{n}-1.
\]
The upper bound is attained for $n=2$ by an odd polygon, and for $n=1$ an
interval has two disjoint facets, so $\chi=1=2^{1}-1$. For $n=3$ the
Four Color Theorem gives $\chi(P^{3})\leq4<7$, so the bound is not
attained (see also \cite{CK26}).

In \cite{Bar25} a simple $4$-polytope with a characteristic map and
$\chi=15$ was constructed, so the bound is attained for $n=4$. The
construction starts from the dual of the cyclic polytope $C^{4}(15)$ and
repairs singular vertices and edges by truncations. It was pointed out in
\cite{Bar25} that the approach does not extend easily to higher dimensions
because of a possible existence of a bad $2$-face, which can be a polygon
with many sides. The question whether the bound is attained for $n\geq5$
was posed as \cite[Problem 3.5]{Bar25}.

In this paper we show that the answer is affirmative for every $n\geq4$.
Throughout, a vertex $v$ is a point of the polytope, whereas $I(v)$ is
the set of facets containing it. Thus, for a simple $n$-polytope,
$|I(v)|=n$ and $\{v\}=\bigcap_{F\in I(v)}F$. We keep this distinction
also for higher-dimensional faces and their defining sets of facets.

\begin{theorem}
\label{Tm_main}For every $n\geq4$ there exists a simple $n$-polytope $P$
admitting a characteristic map over $\mathbb{Z}_{2}$ with
$\chi(P)=2^{n}-1.$ The polytope $P$ can be obtained from
$C^{n}(2^{n}-1)^{\ast}$ by truncations of faces of codimension at least
$3$, with an arbitrary bijection of the facets of $C^{n}(2^{n}-1)^{\ast}$
onto $\mathbb{Z}_{2}^{n}\setminus\{0\}$ as the starting labelling.
\end{theorem}

The main tool is Theorem \ref{Tm_resolution}. It says that an arbitrary
labelling of facets of a simple polytope by vectors of
$\mathbb{Z}_{2}^{n}$ can be turned into a characteristic map by
truncations of faces on which the nonsingularity condition fails. The new
labels may be taken from an arbitrary prescribed spanning set.

A result of this type was proved by Choi and Jang \cite[Lemma 2.5]{CJ26}
in the dual language of stellar subdivisions. They used it in a different
context, for a construction of a polytopal sphere without a mod $3$ characteristic
map, which gives a negative answer to the toric lifting problem. We give a
self-contained proof over an arbitrary field, since we need that only faces
of codimension at least $3$ are truncated, and that the new labels can be
taken from a given set. The latter gives the oriented case. For a proper
starting labelling their technique, applied to the boundary of a cyclic
polytope, already gives the unrestricted statement that the bound $2^{n}-1$
is attained, so we do not claim independence of that step.

Moreover, we construct a simple $n$-polytope with $\chi=2^{n}-1$ which is
the orbit space of both a small cover and a quasitoric manifold, and we obtain
the sharp bound $2^{n-1}$ for oriented small covers for every $n\geq4$. For
$n=4,5,6$ we also report computer-verified examples.

\section{Preliminaries}

All polytopes are considered up to combinatorial equivalence. A polytope
of dimension $n$ is \emph{simple} if each vertex belongs to exactly $n$
facets. We write $I(v)=\{F:v\in F\}$ for this set, so
$\{v\}=\bigcap_{F\in I(v)}F$ and $|I(v)|=n$. Two distinct facets of a simple polytope intersect if and only if
they share a vertex, and then their intersection is a face of codimension
$2$. The \emph{facet graph} of $P$ has the facets as vertices, two facets
being adjacent if they intersect, and $\chi(P)$ denotes its chromatic number.

For a family of facet labels we retain the facet indices: $\nu(S)$ means
the indexed family $(\nu(F))_{F\in S}$, even when some labels coincide.
We use $\bigcap S$ only as shorthand for $\bigcap_{F\in S}F$; it is a
geometric intersection, whereas $S$ is a set of facets.

For a set $S$ of distinct facets with $\bigcap_{F\in S}F\not=\emptyset$, the
face $\bigcap_{F\in S}F$ of a simple polytope $P$ has codimension
$\left\vert S\right\vert $, and the facets of $P$ containing it are
exactly the elements of $S$. The facets of such a face $H$ are the sets
$H\cap F$ for facets $F\not \in S$ which meet $H$, and distinct such $F$
give distinct facets of $H$.

Recall that the cyclic polytope $C^{n}(m)$ is the convex hull of $m$
points on the moment curve. It is simplicial, and for $n\geq4$ every two
of its vertices form an edge, see \cite{Gale63, Ziegler}. Hence for
$n\geq4$ every two facets of the simple polytope $C^{n}(m)^{\ast}$
intersect.

The \emph{truncation} of a face $L=\bigcap_{F\in S}F$ of codimension
$s=\left\vert S\right\vert \geq2$ is the cut by a hyperplane which
separates the vertices of $L$ from the other vertices of $P$. Such a
hyperplane exists since $L$ is a face. The following facts are standard,
see \cite{BP15, ES74}.

The result is a simple polytope $P^{\prime}$ (a convex polytope, since it is
cut by a hyperplane) with the facets of $P$ and
one new facet $F^{\prime}\cong L\times\Delta^{s-1}$. The vertices of $L$
are replaced by vertices $v_{F}$, one for each vertex $v$ of $L$ and
each $F\in S$, determined by
\[
I(v_{F})=(I(v)\setminus\{F\})\cup\{F^{\prime}\}.
\]
These vertices are distinct, so the total number of new vertices is
$s f_{0}(L)$, where $f_{0}(L)$ is the number of vertices of $L$.
All vertices outside $L$ are unchanged.
Dually, this is the stellar subdivision of the corresponding face of the
boundary complex of the simplicial polytope $P^{\ast}$.

\begin{figure}[h]
\begin{center}
\includegraphics[scale=0.9]{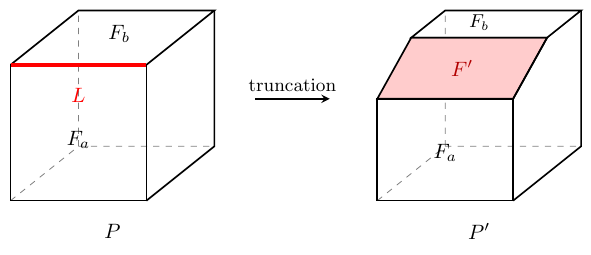}
\end{center}
\caption{The truncation of the edge $L=F_{a}\cap F_{b}$ of a cube $P$ (a face of
codimension $s=2$). The new facet $F^{\prime}\cong L\times\Delta^{1}$ is a
rectangle, the two vertices of $L$ are replaced by the four vertices
with incidence sets $(I(v)\setminus\{F\})\cup\{F^{\prime}\}$,
where $F\in\{F_{a},F_{b}\}$, and all
other vertices are unchanged.}%
\label{Fig_truncation}%
\end{figure}

 For $s=1$, when $L$ is a facet $F$, we use the convention that a
truncation of $L$ is the replacement of $F$ by a new facet $F^{\prime}$ in
the same position. The polytope is unchanged and only the label of the
facet is changed.

\begin{lemma}
\label{Lemma_adjacent}Let $L$ be a face of codimension at least $2$ of a
simple polytope $P$ and let $F_{a},F_{b}$ be distinct intersecting facets.
Their intersection is destroyed by truncating $L$ if and only if
$L=F_{a}\cap F_{b}$. In particular, truncations of codimension at least
$3$ preserve every intersection of two original facets.
\end{lemma}

\begin{proof} The face $F_{a}\cap F_{b}$ has dimension $n-2$. If it is
contained in $L$, the codimension assumption forces $L=F_{a}\cap F_{b}$.
Otherwise it has a vertex outside $L$, which survives the truncation and
still belongs to both facets. If $L=F_{a}\cap F_{b}$, all its vertices
are removed and every new vertex omits one of $F_{a},F_{b}$, so their
intersection is destroyed.
\end{proof}

\section{Resolution of a labelling}

Let $P$ be a simple $d$-polytope and let $U$ be a vector space of dimension
$d$ over a field $F$ (in the main application $F=\mathbb{Z}_{2}$). A \emph{labelling} is an arbitrary map
$\nu$ from the set of facets of $P$ to $U$, where zero and repeated labels
are allowed. A vertex $v$ is \emph{nonsingular} if $\nu(I(v))$ is a basis of
$U$, and $\nu$ is a characteristic map if all vertices are nonsingular.

A \emph{circuit} is a minimal dependent multiset of vectors. A circuit of
size $1$ is a zero label, a circuit of size $2$ is a pair of proportional labels (equal labels if
$F=\mathbb{Z}_{2}$),
and a circuit $c$ spans a space of dimension $\left\vert c\right\vert -1$.
A set $C=\{F_{i_{1}},\ldots,F_{i_{c}}\}$ of distinct facets is called a
\emph{circuit face} if
\[
H=F_{i_{1}}\cap\cdots\cap F_{i_{c}}\not=\emptyset
\quad\hbox{and}\quad
(\nu(F_{i_{1}}),\ldots,\nu(F_{i_{c}}))\text{ is a circuit}.
\]
Thus $C$ denotes the defining set of facets, and $H$ its geometric
intersection, of codimension $c$.

Notice that $\nu$ is a characteristic map if and only if there is no
circuit face. Namely, if $\bigcap C\not =\emptyset$, then every vertex $v$ of the face $\bigcap C$ satisfies
$C\subseteq I(v)$ and is singular. Conversely, the indexed labels at a
singular vertex contain a circuit. Notice also that a circuit face has at
most $d$ elements, so $A=\mathrm{span}(\nu(C))$ is a proper subspace of
$U$ of dimension $\left\vert C\right\vert -1$.

A truncation of a face $L$ with the new label $\mu\in U$ is
\emph{admissible} if $\mu\not \in\mathrm{span}(\nu(I(v)))$ for every
vertex $v$ of $L$. At a new vertex $w$ with
$I(w)=(I(v)\setminus\{F\})\cup\{F^{\prime}\}$, the new facet
$F^{\prime}$ is a \emph{coloop}: its label is not in the span of the
other labels at $w$. In particular, no circuit at $w$ uses $F^{\prime}$.

\begin{lemma}
\label{Lemma_nocircuit}Let $P^{\prime}$ be obtained from $P$ by an
admissible truncation. Then every circuit face of $P^{\prime}$ consists of
old facets and it is a circuit face of $P$. Moreover, a facet created by
an admissible truncation is not contained in any circuit face, also after
any further admissible truncations.
\end{lemma}

\begin{proof} Let $C^{\prime}$ be a circuit face of $P^{\prime}$ and choose
a vertex $w\in\bigcap_{F\in C^{\prime}}F$. If $w$ is unchanged, its
incidence set already contained $C^{\prime}$ in $P$. Otherwise
$I(w)=(I(v)\setminus\{F\})\cup\{F^{\prime}\}$ for a vertex $v$ of
the truncated face. Since $F^{\prime}$ is a coloop at $w$, we have
$C^{\prime}\subseteq I(v)\setminus\{F\}$. Thus
$v\in\bigcap_{G\in C^{\prime}}G$ in $P$, and the same indexed family
of labels was a circuit before the truncation.

For persistence, let $F_{a}$ be a coloop at a vertex $w$. A further
admissible truncation with label $\mu^{\prime}$ either leaves $w$
unchanged or replaces it by vertices with incidence sets
$I(w^{\prime})=(I(w)\setminus\{F_{b}\})\cup\{F^{\prime\prime}\}$.
Only $F_{b}\not=F_{a}$ is relevant when $w^{\prime}\in F_{a}$.
If $F_{a}$ were not a coloop there, then
\[
\nu(F_{a})=a\mu^{\prime}+r,\qquad
r\in\mathrm{span}(\nu(I(w)\setminus\{F_{a},F_{b}\})).
\]
For $a=0$ this contradicts the coloop property at $w$; for $a\not=0$
it implies $\mu^{\prime}=a^{-1}(\nu(F_{a})-r)\in
\mathrm{span}(\nu(I(w)))$, contrary to admissibility. Every vertex of
a surviving facet arises from a vertex of that facet, so the property
persists under all further admissible truncations.
\end{proof}

Let us introduce the notation for a quotient. Let $C$ be a circuit face,
$A=\mathrm{span}(\nu(C))$, $H=\bigcap_{F\in C}F$ and $u\in U\setminus A.$ Then $H$
is a face of dimension $h=d-\left\vert C\right\vert .$ We put
$U_{H}=U/(A+\langle u\rangle)$, which has dimension $h$, and we call the
labelling
\[
\nu_{H}(H\cap F)=\nu(F)+(A+\langle u\rangle)\in U_{H}
\]
of the facets of $H$ the \emph{quotient labelling}. Here $F\not\in C$
is a facet meeting $H$; the value is a coset in $U_{H}$. The face $H$
is determined by the common intersection of the facets in $C$, whereas
its vertices correspond to the incidence sets $I(w)\supseteq C$.
For such a vertex,
$\mathrm{span}(\nu(I(w)))=A+\mathrm{span}(\nu(I(w)\setminus C))$.

Further, if $\overline{\mu}\in U_{H}$ is a label for a truncation in $H$,
then every $\mu\in U$ in the class $\overline{\mu}$ is called a
\emph{lift} of $\overline{\mu}.$

\begin{lemma}
\label{Lemma_quotient}In the above notation the following holds.

\noindent\textup{(i)} The vertices of $P$ contained in $H$ are the
vertices of $H$. If $L$ is a nonempty \emph{proper} face of $H$, then the
truncation of $L$ in $P$ induces the truncation of $L$ in $H$ (for a facet
$L$ of $H$ with the convention for $s=1$), and the set $C$ remains a circuit
face.

\noindent\textup{(ii)} If a truncation of a nonempty proper face
$L\subset H$ with a label
$\overline{\mu}$ is admissible for $\nu_{H}$, then the truncation of $L$
in $P$ with any lift $\mu$ of $\overline{\mu}$ is admissible for $\nu.$

\noindent\textup{(iii)} If $\nu_{H}$ is a characteristic map of $H$, then
the truncation of $H$ itself with the label $u$ is admissible for $\nu.$
\end{lemma}

\begin{proof} (i) Write $L=\bigcap_{F\in S}F$ with $S\supsetneq C$.
At a vertex $v$ of $L$, a new vertex with incidence set
$(I(v)\setminus\{F\})\cup\{F^{\prime}\}$ lies in the surviving face
$H$ if and only if $F\in S\setminus C$. These are exactly the
vertices of the induced truncation in $H$. Since $S\setminus C$ is
nonempty, $H$ survives, and $C$ remains a circuit face.

(ii) If $\mu\in\mathrm{span}(\nu(I(w)))$ for a vertex $w$ of $L$,
its image in $U_{H}$ belongs to the span of the quotient labels of the
facets of $H$ at $w$, since the labels of $C$ have zero image. This
contradicts admissibility of $\overline{\mu}$ in $H$.

(iii) At a vertex $w$ of $H$, the $h$ quotient labels indexed by
$I(w)\setminus C$ form a basis of $U_{H}$. Hence
$A+\langle u\rangle+\mathrm{span}(\nu(I(w)\setminus C))=U$.
If $u\in\mathrm{span}(\nu(I(w)))$, the left-hand side would equal
$\mathrm{span}(\nu(I(w)))$. This is impossible because $\nu(C)$ is
dependent, so the $d$ labels at $w$ span a space of dimension less than $d$.
\end{proof}

A \emph{resolution} of a labelling $\nu$ with labels from $X\subseteq U$
is a finite sequence of admissible truncations such that each truncated
face is contained in $\bigcap C$ (in the current polytope) for a set $C$ of
\emph{original} facets which is a circuit face of $\nu$, all new
labels are from $X$, and the resulting labelling is a characteristic map. If no label of $\nu$ is zero, no truncation of codimension $1$
occurs, so all original labels are kept.

\begin{theorem}
\label{Tm_resolution}Let $P$ be a simple $d$-polytope, $\nu$ a labelling
of its facets by vectors of $U$ and $X\subseteq U$ a set which spans $U$.
Then $\nu$ has a resolution with labels from $X.$
\end{theorem}

\begin{proof} We use induction on $d$ and, for a fixed $d$, induction on
the number $c$ of circuit faces. For $d=0$ there is nothing to prove, and
for $c=0$ the labelling is already a characteristic map. So let $c\geq1.$
Let us first fix the data. Let $C$ be a circuit face,
$A=\mathrm{span}(\nu(C))$ and $H=\bigcap C.$ Since $A$ is a proper
subspace and $X$ spans $U$, there is $u\in X\setminus A.$ Let $\nu_{H}$ be
the quotient labelling of $H.$ The image $X_{H}$ of $X$ in $U_{H}$ spans
$U_{H}$, and $\mathrm{dim}(H)<d.$

Let us next resolve the labelling on $H.$ By the induction hypothesis
applied to $H$ with the labelling $\nu_{H}$ and the set $X_{H}$ there is a sequence of admissible
truncations of proper faces of $H$ with labels from $X_{H}$ (empty if $H$ is
a point), after which
$\nu_{H}$ is a characteristic map. By Lemma \ref{Lemma_quotient} (i) and
(ii), each of these truncations lifts to an admissible truncation in $P$
of a face contained in $H$ (a codimension-$1$ step in $H$ lifts to a
truncation of codimension $|C|+1$ in $P$). The lift of the label can be chosen in $X$,
since $X_{H}$ is the image of $X$. After the last of them $\nu_{H}$ is
nonsingular on the current $H.$

Now we truncate $H$ itself. By Lemma \ref{Lemma_quotient} (iii) the
truncation of $H$ with the label $u\in X$ is admissible. Notice that after
this truncation no vertex has an incidence set containing $C$. Namely, the vertices of the old
$H$ are removed, and the new vertices
have incidence sets $(I(v)\setminus\{F\})\cup\{F^{\prime}\}$ with
$F\in C$, none of which contains $C$.

Let us now count the circuit faces. By Lemma \ref{Lemma_nocircuit} none of
the truncations created a circuit face, and every circuit face which
exists now was a circuit face before. Moreover, it is different from $C.$
So the number of circuit faces is smaller than $c.$

Finally, by the induction on $c$ the labelling can be completed to a
characteristic map by further truncations of faces contained in
intersections $\bigcap_{F\in C}F$ defined by circuit faces $C$ of the
current polytope. These sets $C$ are circuit faces of the original
polytope by Lemma \ref{Lemma_nocircuit}, and all the labels are from $X.$
\end{proof}

\begin{proposition}
\label{Prop_survive}Let $\nu$ have no zero labels and let $P$ be obtained from
$P_{0}$ by a resolution of $\nu$. Let $S$ be a set of original facets with
$\bigcap S\not =\emptyset$ in $P_{0}$. Then $\bigcap S\not =\emptyset$ in $P$
if and only if $\nu(S)$ is linearly independent.
\end{proposition}

\begin{proof} Original facets are kept with their labels. If
$\bigcap S\not =\emptyset$ in $P$, then $S\subseteq I(v)$ for a vertex $v$ of $P$, whose
labels form a basis, so $\nu(S)$ is independent. Conversely let $\nu(S)$ be
independent and let $L=\bigcap S_{L}$ be a truncated face, where
$S_{L}\supseteq C$ for an original circuit face $C$. Then $S_{L}\not \subseteq S$,
since $S$ contains no circuit. Choose a vertex $v$ of the current polytope with $S\subseteq I(v)$.
If $v\not\in L$, it survives. Otherwise $S_{L}\subseteq I(v)$, and
for $F\in S_{L}\setminus S$ a new vertex $w$ has
$I(w)=(I(v)\setminus\{F\})\cup\{F^{\prime}\}\supseteq S$. So $\bigcap S$ stays
non-empty after each truncation.
\end{proof}

The circuit faces counted in the proof of Theorem \ref{Tm_resolution}
are sets of actual facets, not sets of vector values. Distinct facet sets
can carry the same indexed vector values, and are counted separately.
Lemma \ref{Lemma_nocircuit} excludes new circuit faces even when labels
repeat. No global injectivity of the labelling is used.

For $|S|=2$ this is Lemma \ref{Lemma_adjacent} for resolutions: two
intersecting original facets keep intersecting if and only if their labels are
not proportional.

\noindent\textbf{Proof of Theorem \ref{Tm_main}.} Let $n\geq4$,
$m=2^{n}-1$, $P_{0}=C^{n}(m)^{\ast}$ and $U=\mathbb{Z}_{2}^{n}.$ Let $\nu$
be an arbitrary bijection of the facets of $P_{0}$ onto $U\setminus\{0\}.$
All labels are non-zero and pairwise different, so every circuit has at
least $3$ elements. Hence the intersection defined by a circuit face has
codimension at least $3$.

We apply Theorem \ref{Tm_resolution} with $X=U\setminus\{0\}$ and denote
by $P$ the obtained polytope. Every truncated face is contained in
$\bigcap_{F\in C}F$ for a circuit face $C$ of $P_{0}$, so its codimension
is at least $3$. By Lemma
\ref{Lemma_adjacent} the $m$ original facets, which pairwise intersect in
$P_{0}$ since $C^{n}(m)$ is $2$-neighborly, pairwise intersect in $P.$

These facets need pairwise different colours, so $\chi(P)\geq2^{n}-1.$ The
final labelling is a characteristic map, hence a proper colouring with at
most $2^{n}-1$ colours, and so $\chi(P)=2^{n}-1.$ \hfill $\square$

The proof did not use $F=\mathbb{Z}_{2}$, so the same argument works over
any field. Let $q$ be a prime power. A characteristic map over
$\mathbb{F}_{q}$ is a labelling with a basis of $\mathbb{F}_{q}^{n}$ at every
vertex; here we use it as a purely algebraic notion. Two facets sharing a
vertex have non-proportional labels, so the facets are properly coloured by
the points of the projective space, and $\chi\leq N_{q}(n)=(q^{n}-1)/(q-1).$

\begin{corollary}
\label{Cor_Fq}For every prime power $q$ and $n\geq4$ the largest chromatic
number of a simple $n$-polytope with a characteristic map over
$\mathbb{F}_{q}$ is $N_{q}(n)$, and every integer between $n$ and $N_{q}(n)$
is attained.
\end{corollary}

\begin{proof} Take $X$ with one representative of each projective point and
a bijection of the facets of $C^{n}(N_{q}(n))^{\ast}$ onto $X$. No label is
zero and no two are proportional, so every circuit has at least $3$
elements, and the proof of Theorem \ref{Tm_main} applies verbatim. For
$n+1\leq k\leq N_{q}(n)$ take $k$ representatives containing a basis and
$C^{n}(k)^{\ast}$ as the starting polytope; for $k=n$ take the cube with
labels $e_{1},\dots,e_{n}$ on the pairs of opposite facets.
\end{proof}

In particular, for every $n\geq4$ and prime powers $q>r$, there is a
simple $n$-polytope with a characteristic map over $\mathbb{F}_{q}$ and
none over $\mathbb{F}_{r}$, since $N_{q}(n)>N_{r}(n)$.

This has a formulation in terms of the universal complexes studied by
Barali\'{c}, Vavpeti\v{c} and Vu\v{c}i\'{c} \cite{BVV23,BVV24}.
Write $X(\mathbb{F}_{q}^{\ell})$ for the simplicial complex whose vertices
are the non-zero vectors and whose simplices are the linearly independent
sets. A simplicial map is \emph{nondegenerate} if it is injective on each
simplex. A characteristic map on $P$ is equivalent to a nondegenerate
simplicial map $K=\partial P^{\ast}\longrightarrow X(\mathbb{F}_{q}^{n})$.
In this language Theorem \ref{Tm_resolution} extends a vertex labelling
to such a map after stellar subdivisions, with prescribed new labels.

For a prime $p$, the mod $p$ Buchstaber invariant of a complex $K$ on
$m$ vertices is $s_{p}(K)=m-\ell_{p}(K)$, where $\ell_{p}(K)$ is the least
$\ell$ admitting a nondegenerate map to $X(\mathbb{F}_{p}^{\ell})$
\cite{BVV24}. Our construction gives, for every pair of primes $p>r$
and every $n\geq4$, a polytopal $(n-1)$-sphere $K$ for which
\[
s_{p}(K)=m-n> s_{r}(K).
\]
Indeed, take $P$ from Corollary \ref{Cor_Fq} with $q=p$ and
$\chi(P)=N_{p}(n)$. Its dual sphere maps nondegenerately to
$X(\mathbb{F}_{p}^{n})$, and its dimension forces $\ell_{p}(K)=n$.
A map to $X(\mathbb{F}_{r}^{n})$ would colour the facets by at most
$N_{r}(n)$ projective points, a contradiction. Thus
$\ell_{r}(K)\geq n+1$, proving the strict inequality. This yields the
separation on polytopal spheres in every dimension at least $3$.

The following refinement does not need a complete starting graph.

\begin{theorem}
\label{Tm_Q}Let $Q$ be a simple $n$-polytope with $\chi(Q)=k$, and let $q$ be a
prime power with $k\leq N_{q}(n)$. Then there is a simple $n$-polytope $P$,
obtained from $Q$ by truncations of faces of codimension at least $3$, with a
characteristic map over $\mathbb{F}_{q}$ and $\chi(P)=k$. The original facets
induce in $P$ the same intersection graph as in $Q$.
\end{theorem}

\begin{proof} Since $n$ facets meet at a vertex, $k\geq n$. Let $X$ be a set
of $k$ projective representatives containing a basis, and label the facets of
$Q$ by $X$ according to a proper colouring with $k$ colours. There are no zero
labels, and a circuit of size $2$ would be two intersecting facets with
proportional, hence equal, labels, which the colouring excludes. So every
circuit face has at least $3$ facets. Theorem \ref{Tm_resolution} with $X$
truncates only faces of codimension at least $3$, and by Lemma
\ref{Lemma_adjacent} no intersection of original facets is lost. A truncation
creates no intersection between two formerly disjoint old facets, since a new
vertex consists of a subset of an old vertex together with the new facet.
Hence the old facets induce the intersection graph of $Q$, so $\chi(P)\geq
k$, and all labels lie in $X$, so $\chi(P)\leq k$.
\end{proof}

For $q=2$ this says that every simple polytope with $\chi\leq2^{n}-1$ has a
small-cover refinement with the same chromatic number, and a maximal clique is
not needed for the lower bound.

\begin{remark}
In the proof of Theorem \ref{Tm_resolution} the labels of new facets are
chosen only from the condition $u\in X\setminus A$ for one proper subspace
$A$, and its lifts. This is why they can be taken from any spanning set.
In the dual language, for $F=\mathbb{Z}_{2}$ and $X=U$ Theorem
\ref{Tm_resolution} is a form of \cite[Lemma 2.5]{CJ26}.
\end{remark}

\section{The quasitoric case}

A non-zero integer vector is \emph{primitive} if the greatest common
divisor of its coordinates is $1$, equivalently if it is not an integer
multiple $k u$ with $k>1$ and $u\in\mathbb{Z}^{n}$.
A characteristic map over $\mathbb{Z}$ assigns primitive vectors to
facets so that the vectors at every vertex form a lattice basis, or
equivalently have determinant $\pm1$. Over a field, linear independence
of $n$ vectors already gives a basis. Over $\mathbb{Z}$, a non-zero
determinant gives only a sublattice of finite index. We therefore first
apply the field resolution over $\mathbb{Q}$ and then reduce these
indices to $1$. It gives a quasitoric manifold, and its reduction
modulo $2$ is a characteristic map over $\mathbb{Z}_{2}$. Notice that a
given $\mathbb{Z}_{2}$ characteristic map does not always lift to
$\mathbb{Z}$ \cite{CJ26}. Here we construct the integral map together with
the polytope; reducing a determinant by subdividing with a lattice point of
the fundamental parallelepiped is a standard device of toric resolution, see
\cite[\S2.6]{Ful93} and \cite{BvT17}, and what is needed here is that the intersections of old facets
are protected; truncations of codimension $2$ are allowed, but no
intersection of two old facets is removed. The polytope obtained is in
general different from the one in Theorem \ref{Tm_main}.

\begin{lemma}
\label{Lemma_saturated}Let $a\not =b$ be non-zero vectors of
$\{0,1\}^{n}\subseteq\mathbb{Z}^{n}.$ Then $a,b$ can be extended to a
basis of $\mathbb{Z}^{n}$. In particular,
$(\mathbb{Q}a+\mathbb{Q}b)\cap\mathbb{Z}^{n}=\mathbb{Z}a+\mathbb{Z}b.$
\end{lemma}

\begin{proof}
Two different non-zero $0/1$ vectors are not parallel, so some $2\times2$
minor of the matrix $(a,b)$ is non-zero. Each of them is in $\{-1,0,1\}$,
so the gcd of the $2\times2$ minors is $1.$ This is the criterion for a
pair to be extendable to a basis of $\mathbb{Z}^{n}$, by the Smith normal
form.
\end{proof}

\begin{theorem}
\label{Tm_Z}Let $P_{0}$ be a simple $n$-polytope and let $V$ assign to the
facets of $P_{0}$ non-zero vectors of $\{0,1\}^{n}$ such that
intersecting facets get distinct vectors. Then there is a polytope $P$ obtained
from $P_{0}$ by truncations of faces of codimension at least $2$, and a
characteristic map $\Lambda$ over $\mathbb{Z}$ on $P$ extending $V$, such
that two facets of $P_{0}$ intersect in $P$ if and only if they intersect in
$P_{0}$. In particular,
for every $n\geq4$ there is a simple $n$-polytope $P$ with
$\chi(P)=2^{n}-1$ which is the orbit space of a quasitoric manifold, and
of a small cover.
\end{theorem}

\begin{proof} The facets of $P_{0}$ are called \emph{original}. We keep their
labels, denote the current labelling by $V$, and write $V(S)$ for the
family of column vectors indexed by a set $S$ of facets. All labels
chosen below are primitive integer vectors. We follow the field
resolution and then perform the extra lattice step.

\medskip
\noindent\textbf{Claim A.} \emph{Truncations of codimension at least $3$
give a labelling with non-zero determinant at every vertex.}\smallskip

\noindent Apply Theorem \ref{Tm_resolution} over $\mathbb{Q}$ with the
spanning set $X$ of all primitive vectors of $\mathbb{Z}^{n}$. The labels
of intersecting original facets are distinct non-zero $0/1$ vectors,
hence are not proportional over $\mathbb{Q}$. Thus every original
$\mathbb{Q}$-circuit face $C$ has at least three facets. Each face
truncated in the resolution lies in the geometric face
$\bigcap_{F\in C}F$ for such a $C$, so has codimension at least $3$.
The resulting labels at every vertex form a basis over $\mathbb{Q}$,
and consequently have non-zero integer determinant. All original labels
are preserved. This establishes Claim A.

\medskip
\noindent\textbf{Claim B.} \emph{Further truncations make all determinants
equal to $\pm1$, preserving the intersections of original facets.}\smallskip

\noindent For a vertex $v$ put
$D(v)=|\det V(I(v))|$. By Claim A these are positive integers. If some
$D(v)>1$, choose $v$ with maximal $D(v)$ and write
$I(v)=\{F_{1},\ldots,F_{n}\}$. The sublattice
$L_{v}=\sum_{j=1}^{n}\mathbb{Z}V(F_{j})$ has index $D(v)$ in
$\mathbb{Z}^{n}$. Its half-open fundamental parallelepiped is
\[
\Pi_{v}=\left\{\sum_{j=1}^{n}c_{j}V(F_{j}):0\leq c_{j}<1\right\}.
\]
We use the following standard fact (see, e.g., \cite{BvT17}); we include
the short proof. \emph{Every coset of $L_{v}$ in $\mathbb{Z}^{n}$ has exactly
one representative in $\Pi_{v}$.} Indeed, let $x\in\mathbb{Z}^{n}$ and write
$x=\sum_{j}a_{j}V(F_{j})$; the $a_{j}$ are uniquely determined and rational,
since $V(I(v))$ is an invertible integer matrix. Then
$x^{\ast}=x-\sum_{j}\lfloor a_{j}\rfloor V(F_{j})$ lies in $x+L_{v}$ and in
$\Pi_{v}$, its coefficients being the fractional parts of the $a_{j}$. If
$x^{\ast},x^{\ast\ast}\in\Pi_{v}$ are in the same coset, their difference lies
in $L_{v}$, so it has integer coefficients, and these lie in $(-1,1)$; hence
the difference is $0$. Since $[\mathbb{Z}^{n}:L_{v}]=D(v)>1$, there is a
coset different from $L_{v}$, and its representative $w_{0}\in\Pi_{v}$ is an
integer vector which is non-zero (the coset of $0$ is $L_{v}$ itself).

Let $k\geq1$ be the greatest common divisor of the coordinates of $w_{0}$ and
put $w=w_{0}/k$. Then $w\in\mathbb{Z}^{n}$ is non-zero and primitive, by the
definition of primitive vectors above. Its coefficients with respect to
$V(F_{1}),\ldots,V(F_{n})$ are $c_{j}=c_{j}^{0}/k\in[0,1)$, so $w\in\Pi_{v}$ as
well. (We do not need that $w$ is in the coset of $w_{0}$; we only use that
$w$ is a non-zero primitive integer vector in $\Pi_{v}$.)

Let $T=\{F_{j}:c_{j}>0\}$, a set of facets, and let
$H=\bigcap_{F\in T}F$, a nonempty geometric face containing $v$.
The set $T$ has at least two elements: otherwise an integer vector would
be a multiple $c_{j}V(F_{j})$ of a primitive vector with $0<c_{j}<1$.
This is impossible, since the coordinates of a primitive vector have
an integer linear combination equal to $1$.
Moreover, $T$ cannot consist of two original facets $F_{a},F_{b}$.
They intersect at $v$, and Lemma \ref{Lemma_saturated} gives
\[
(\mathbb{Q}V(F_{a})+\mathbb{Q}V(F_{b}))\cap\mathbb{Z}^{n}
=\mathbb{Z}V(F_{a})+\mathbb{Z}V(F_{b}).
\]
Uniqueness of the rational coefficients would force $c_{a},c_{b}$ to
be integers, contrary to $0<c_{a},c_{b}<1$.

Truncate $H$ and label the new facet $F^{\prime}$ by $w$. For a vertex
$z$ of $H$ and $F_{i}\in T$, the corresponding new vertex $z_{i}$ has
$I(z_{i})=(I(z)\setminus\{F_{i}\})\cup\{F^{\prime}\}$. Let $M_{z}$ be the
$n\times n$ matrix whose columns are $V(F)$, $F\in I(z)$, in a fixed order,
and let $M_{z_{i}}$ be obtained from $M_{z}$ by replacing the column
$V(F_{i})$ with $w$. Recall that
$w=\sum_{F_{j}\in T}c_{j}V(F_{j})$, where the coefficients $c_{j}$ were
determined at the vertex $v$ and do not depend on $z$. Since
$T\subseteq I(z)$, every column $V(F_{j})$ with $F_{j}\in T$ is a column of
$M_{z}$. By the linearity of the determinant in the replaced column,
\[
\det M_{z_{i}}=\sum_{F_{j}\in T}c_{j}\det M_{z}^{(i\to j)},
\]
where $M_{z}^{(i\to j)}$ is $M_{z}$ with the column $V(F_{i})$ replaced by
$V(F_{j})$. For $j\neq i$ the matrix $M_{z}^{(i\to j)}$ has two equal
columns (the column $V(F_{j})$ already occurs in $M_{z}$), so its
determinant is $0$, while for $j=i$ it is $M_{z}$. Hence
\[
\det M_{z_{i}}=c_{i}\det M_{z},\qquad\text{so}\qquad
D(z_{i})=c_{i}D(z),\quad 0<D(z_{i})<D(z),
\]
because $0<c_{i}<1$ for $F_{i}\in T$ and $D(z)\geq1$. The new determinant
is an integer, as all columns are integer vectors, and $D(z_{i})\geq1$.
All vertices outside $H$ are unchanged.

This procedure terminates. At a step with maximal determinant $M>1$,
at least one vertex with determinant $M$ disappears and no such vertex
is created. After finitely many steps with that maximum it decreases.
The maximum is a positive integer, so eventually all determinants are $1$.

By Lemma \ref{Lemma_adjacent}, a truncation can destroy an intersection
of two original facets only by truncating exactly that intersection.
Claim A uses codimension at least $3$, while in Claim B this would mean
that $T$ is a pair of original facets, which was excluded. Truncations
also create no intersection of formerly disjoint facets. This
establishes Claim B and the first statement of the theorem.

\medskip
For the second statement take $P_{0}=C^{n}(2^{n}-1)^{\ast}$ with $n\geq4$
and a bijection of its facets onto the non-zero vectors of $\{0,1\}^{n}$.
Its $2^{n}-1$ original facets still pairwise intersect, so
$\chi(P)\geq2^{n}-1$. The reduction of the integral characteristic map
modulo $2$ is characteristic, since every determinant is $\pm1$, and
therefore $\chi(P)\leq2^{n}-1$.
\end{proof}

\begin{corollary}
\label{Cor_Zrel}A simple $n$-polytope $Q$ has a refinement $P$ with a
characteristic map over $\mathbb{Z}$, in which two original facets intersect
if and only if they intersect in $Q$, if and only if $\chi(Q)\leq2^{n}-1$.
\end{corollary}

\begin{proof} If $\chi(Q)\leq2^{n}-1$, assign distinct non-zero vectors of
$\{0,1\}^{n}$ to the colour classes of a proper colouring and apply Theorem
\ref{Tm_Z}. Conversely, the reduction modulo $2$ of an integral characteristic
map on $P$ properly colours $P$, hence the original facets, with at most
$2^{n}-1$ colours.
\end{proof}

This does not say that $\chi(P)=\chi(Q)$, since the new integral labels may
need further colours.

For $n=5$ the lift can be found also directly. For the polytope with $174$
facets and $1636$ vertices there is a map $\Lambda$ with
$\Lambda\equiv\lambda\pmod 2$ and determinant $\pm1$ at all vertices. We
assign one integer vector to each of the $31$ non-zero vectors of
$\mathbb{Z}_{2}^{5}$, then $1358$ different bases occur at the vertices,
and a search over the signs of entries gives a solution. We verified it by
exact integer arithmetic.

\section{The oriented case}

A small cover is orientable if and only if, for a suitable choice of a
basis of $\mathbb{Z}_{2}^{n}$, all values of the characteristic map have
an odd number of non-zero coordinates \cite{NN05}. The vectors of odd
weight span $\mathbb{Z}_{2}^{n}$, since they contain the basis vectors,
and there are $2^{n-1}$ of them. So in the oriented case $\chi\leq2^{n-1}.$

\begin{proposition}
\label{Prop_or}For every $n\geq4$ there is a simple $n$-polytope $P$ with
$\chi(P)=2^{n-1}$ which is the orbit space of an oriented small cover. Every
integer between $n$ and $2^{n-1}$ occurs as such a chromatic number.
\end{proposition}

\begin{proof} Let $X$ be the set of vectors of $\mathbb{Z}_{2}^{n}$ of odd
weight, so $\left\vert X\right\vert =2^{n-1}$, and let
$P_{0}=C^{n}(2^{n-1})^{\ast}$ with an arbitrary bijection $\nu$ of its
facets onto $X.$ We apply Theorem \ref{Tm_resolution} with this $X.$

All labels are non-zero and pairwise different, so circuit faces have at
least $3$ elements. Moreover, the sum of an odd number of vectors of odd
weight has odd weight, so it is non-zero, and the circuits have an even
number of elements, hence at least $4$, so all truncations have codimension at
least $4$. By the same argument as in the proof of Theorem
\ref{Tm_main} the $2^{n-1}$ old facets pairwise intersect in the final
polytope $P$, so $\chi(P)\geq2^{n-1}.$

All labels of the characteristic map are from $X$, so the small cover is
orientable by \cite{NN05}, and $\chi(P)\leq\left\vert X\right\vert%
=2^{n-1}.$
\end{proof}

For $n=4$ and $n=5$ this reproves \cite[Theorems 3.6, 3.7]{Bar25}, and the
cases $n\geq6$ are new. Computer-verified examples exist for $n=4,5,6,7$; for $n=7$ the polytope has
$13391$ facets and $340374$ vertices, and $\chi=64.$

The same argument gives the full spectrum. For $n\leq k\leq2^{n-1}$ take $k$
vectors of odd weight which contain the standard basis; if $k\geq n+1$ use the
starting polytope $C^{n}(k)^{\ast}$, and if $k=n$ use the cube with labels
$e_{1},\dots,e_{n}$ on opposite pairs of facets. The resulting oriented small
cover has $\chi=k$. In the unrestricted case, with an arbitrary spanning set
$X\subseteq\mathbb{Z}_{2}^{n}\setminus\{0\}$ with $\left\vert X\right\vert
\geq n+1$ and the starting polytope $C^{n}(\left\vert X\right\vert )^{\ast}$,
we get a characteristic map with image $X$ and $\chi=\left\vert X\right\vert$;
so every value of $\chi$ between $n$ and $2^{n}-1$ is attained for $n\geq4$
(the value $n$ by the cube).

More generally, let $Q$ be any simple $n$-polytope with $\chi(Q)=k\leq2^{n-1}$.
Label its facets according to a proper colouring with $k$ colours by $k$
vectors of odd weight which contain the standard basis. No two intersecting
facets have equal labels, and the sum of an odd number of odd-weight vectors is
non-zero, so every circuit face has an even number of facets, hence at least
$4$. Theorem \ref{Tm_resolution} gives a refinement $P$ of $Q$ by truncations of
codimension at least $4$, with an orientable small cover and $\chi(P)=k$, in
which the original facets induce the intersection graph of $Q$, as in Theorem
\ref{Tm_Q}.

\section{Concluding remarks}

We showed that the upper bound $\chi(P^{n})\leq2^{n}-1$ of Davis and
Januszkiewicz is attained for every $n\geq4$, which answers \cite[Problem
3.5]{Bar25}. Together with $n=1,2$ and the failure for $n=3$ this settles
the question in all dimensions. The polytopes are obtained from
$C^{n}(2^{n}-1)^{\ast}$ by truncations of faces of codimension at least
$3$. Moreover, there is a polytope with $\chi=2^{n}-1$ which is the orbit space of
both a small cover and a quasitoric manifold (Theorem \ref{Tm_Z}), and in the
oriented case the sharp bound $2^{n-1}$ is attained for every $n\geq4$.

A labelling of the facets of $C^{n}(2^{n}-1)^{\ast}$ can have singular
vertices of corank at least $2$, so a repair of single vertices is not
enough, and this is why circuit faces in Theorem \ref{Tm_resolution} are
repaired recursively.

\begin{figure}[h]
\begin{center}
\includegraphics[scale=0.9]{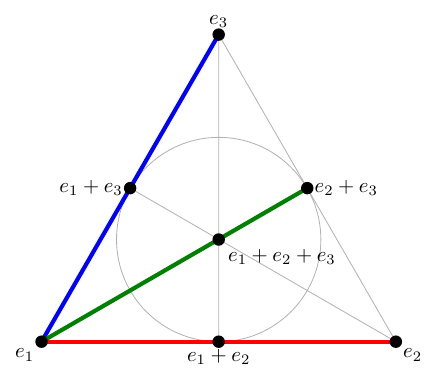}
\end{center}
\caption{The Fano plane $\mathbb{Z}_{2}^{3}\setminus\{0\}$. For a hexagon whose
sides have the classes $e_{1},e_{2},e_{1},e_{3},e_{1},e_{2}+e_{3}$ the lines
$\mathrm{span}(\mu_{j},\mu_{j+1})$ are the three bold lines through $e_{1}$.
They cover all seven points, so no label $u$ is admissible for a direct
truncation of the hexagon.}%
\label{Fig_Fano}%
\end{figure}

Let us explain why directly truncating a bad face may fail. Let $G=F_{a}\cap F_{b}\cap F_{c}$ be a bad $2$-face in dimension $5$ whose
labels form a circuit, so $W=\mathrm{span}\{\nu(F_{a}),\nu(F_{b}),\nu(F_{c})\}$
has dimension $2$, suppose that $G$ is a polygon, and let
$\mu_{1},\ldots,\mu_{k}$ be the classes of the labels of its sides in
$\mathbb{Z}_{2}^{5}/W\cong\mathbb{Z}_{2}^{3}$, i.e., points of the Fano
plane, with consecutive classes distinct and non-zero. A truncation of $G$ with a label $u$ is admissible if and only if
the class of $u$ lies on none of the lines
$\mathrm{span}(\mu_{j},\mu_{j+1})$. These lines can cover the whole Fano
plane, as in Figure \ref{Fig_Fano}, so a direct truncation is impossible.
This is the obstruction pointed out in \cite{Bar25}, and the recursion in
Theorem \ref{Tm_resolution} repairs it by truncating faces inside of $G$
first.

We also note that a bad edge $F_{3}\cap F_{6}\cap F_{7}$ is missing in
\cite[Example 3.2]{Bar25}, since the labels $e_{4},$ $e_{2}+e_{4}$ and
$e_{2}$ sum to zero. So the example has $14$ bad edges and the face vector
of the resulting polytope is $(197,394,243,46)$, instead of
$(193,386,228,45)$ which does not satisfy the Euler relation. This does
not influence Theorem 3.1 of \cite{Bar25}.

The construction was also implemented on a computer by several independent
programs, and the examples below were obtained in this way and checked
by replaying the truncations. With the
cyclic starting polytope, the strongly admissible recursion of Theorem
\ref{Tm_resolution} gives examples with $38$ facets for $n=4$, with $164$
facets for $n=5$ (face vector $(1636,4090,3596,1304,164)$), and with $6681$
facets for $n=6$. A more flexible repair rule gives a polytope with $146$
facets for $n=5$, with face vector $(1578,3945,3444,1221,146)$. For $n=4$ a
search starting from a non-cyclic neighborly $4$-polytope with $15$ vertices
gives a simple polytope with $20$ facets, $110$ vertices, $\chi=15$ and a
characteristic map whose $0/1$ lift is integral (all $110$ determinants are
$\pm1$), so it is the orbit space of both a small cover and a quasitoric
manifold. These numbers are upper bounds and not known minima. The examples
were found by computer, and we do not include certificates for them.
\begin{remark}
\label{Rem_f4}For $n\geq4$, let $f_{\mathbb{Z}_{2}}(n)$ and $f_{\mathbb{Z}}(n)$ be the
smallest numbers of facets of a simple $n$-polytope with $\chi=2^{n}-1$ and a
characteristic map over $\mathbb{Z}_{2}$, respectively over $\mathbb{Z}$.
These functions count all facets, including those introduced by
truncation, rather than circuit faces or different label values.
Reduction modulo $2$ gives $f_{\mathbb{Z}_{2}}(n)\leq f_{\mathbb{Z}}(n)$. For
$n=4$ we have $15\leq f_{\mathbb{Z}_{2}}(4)=f_{\mathbb{Z}}(4)\leq20$. Namely,
by the $20$-facet example a minimizer for $\mathbb{Z}_{2}$ has at most $20$
facets, so its dual boundary is a simplicial $3$-sphere with at most $20$
vertices, and every mod $2$ characteristic map on such a sphere lifts to
$\mathbb{Z}$ by \cite[Theorem 3.1]{CGJL26}.
\end{remark}

There is a necessary condition for an example with exactly
$N=2^{n}-1$ pairwise intersecting facets. Every mod $2$ characteristic
map must label them bijectively by the non-zero vectors. Thus, in its
dual boundary, each triple with labels $a,b,a+b$ is a missing triangle:
all its edges are present, but the triangle is absent by linear
dependence. There are $N(N-1)/6$ such triples, since every pair of
distinct non-zero vectors determines one of them, and each triple is
counted by three pairs. This restricts the $2$-neighborly candidates
without requiring them to be neighborly; it does not settle existence.

\bigskip

\medskip\noindent\textbf{AI declaration.}~~Anthropic Claude.ai was used in the
preparation of this manuscript. It assisted with the arguments and their
exposition, and with the verification of the proofs. All results were verified
by the authors, who take full responsibility for the content.

\medskip\noindent\textbf{Acknowledgements.}~~Djordje Barali\'{c} was partially supported by the Serbian Ministry of Science, Technological Development, and Innovation through the Mathematical Institute of the Serbian Academy of Sciences and Arts. ~~R.~\v{S}krekovski acknowledges
partial support from the Slovenian Research and Innovation Agency program
P1--0383.

\end{document}